\documentclass{amsart}

\usepackage{amscd,amsmath}

\newcommand{\f}{\varphi}

\newcommand{\e}{\hat\eta}

\newcommand{\ov}{\overline}

\newcommand{\Ll}{\mathcal{L}}

\newcommand{\CC}{\mathbb{C}}   

\newcommand{\HH}{\mathbb{H}}   

\newcommand{\RR}{\mathbb{R}}   

\newcommand{\ZZ}{\mathbb{Z}}

\newcommand{\HP}[1]{\mathbb{H}P^{#1}}

\numberwithin{equation}{section}

\newtheorem{te}{Theorem}[section]

\newtheorem{pr}{Proposition}[section]

\newtheorem{co}{Corollary}[section]

\theoremstyle{definition}

\theoremstyle{remark}

\newtheorem{re}{Remark}[section]

\begin{document}
\title
{Complex structures on some Stiefel manifolds}
\author{Liviu Ornea}
\address{University of Bucharest, Faculty of Mathematics,\newline 
14
Academiei str. 70109 Bucharest, Romania}
\email{lornea@imar.ro}

\author{Paolo Piccinni}
\address{Universit\`{a} degli Studi di Roma ''La Sapienza''\newline 
 Piazzale Aldo Moro 2, I-00185 Roma, Italy  }
\email{piccinni@mat.uniroma1.it}

 %\date{\today}
\dedicatory{Dedicated to the memory of Prof. Gheorghe Vr\u anceanu}

\begin{abstract} We give some applications of a construction, appeared in
\cite{orp}, of an integrable complex structure on the total
space of an induced Hopf $S^3$-bundle over a Sasakian manifold. 
We show how this construction allows to obtain an uncountable family of inequivalent 
complex structures on the Stiefel manifolds
$V_2({\CC}^{n+1})$ and  $\widetilde V_4({\RR}^{n+1})$, as well as on
some special Stiefel manifolds related to the groups $G_2$ 
and $Spin(7)$. In the case of $V_2({\CC}^{n+1})$, these complex
structures are not compatible
with its standard hypercomplex structure. 
\end{abstract} 
\maketitle
\section{Introduction}
The Stiefel manifolds $V_k({\CC}^{n+1})$ of orthonormal $k$-frames in 
${\CC}^{n+1}$ 
and $\widetilde V_{2h}({\RR}^{n+1})$ of oriented orthonormal real 
$2h$-frames in ${\RR}^{n+1}$ appear in the classical work of
H.-C. Wang \cite{wa} as examples of compact manifolds that admit an infinite 
family of inequivalent homogeneous complex
structures, described by a real parameter. 
The class of homogeneous manifolds with this property, 
all non-k\"ahlerian, includes
also compact simple Lie groups, studied in this respect almost 
simultaneously by H. Samelson \cite{Sam}. More recently, both
Wang's and Samelson's work inspired remarkable constructions of 
both homogeneous and inhomogeneous hypercomplex structures on
some classes of manifolds that 
include $V_2({\CC}^{n+1})$: \cite{Jo}, \cite{Bo-Ga-Ma 2}.  

The aim of this note is to present a simple construction of an 
infinite family of homogeneous complex
structures - this time described by a complex parameter - on 
$V_2({\CC}^{n+1})$, $\widetilde V_4({\RR}^{n+1})$, $G_2$ and
$Spin(7)/Sp(1)$, the latter two manifolds being special 
Stiefel with respect to the geometry of the Cayley numbers. This
construction was suggested by our work on the geometry of the zero 
level set of some moment maps defined on the
quaternionic projective space
\cite{orp}, where some diagrams involve these Stiefel 
manifolds and a definition on them of a complex structure in
the Calabi-Eckmann spirit is quite natural. 

The construction of a complex structure on the total space of an induced 
Hopf $S^1$-bundle over a Sasakian manifold is classical. Locally, one makes the
product between the Sasakian structure of the base manifold and the standard one
of the circle. By a similar technique, one can produce a hypercomplex structure on the total space of a framed
$S^1$-bundle  over a 3-Sasakian manifold, and this is the way to obtain the uncountably many
hypercomplex structures on $V_2({\CC}^{n+1})$ \cite{Bo-Ga-Ma 2}. 

In  Proposition \ref{unu} below we carry on a similar construction,
but using one of the Sasakian structures in the the standard 
$3$-Sasakian structure of the unit sphere $S^3$ that is the fiber of an 
induced Hopf bundle over a Sasakian manifold.

This construction can be applied to some Sasakian manifolds that turn out to be the zero level sets of some
moment maps defined on the quaternionic projective space. Again, the Sasakian 
structures are here induced by that of a sphere,
now $S^{2N+1}$, by means of appropriate induced Hopf circle bundles.  
With these Sasakian manifolds as base spaces, our
method provides a complex structure on the Stiefel manifolds
$V_2({\CC}^{n+1})$, $\widetilde V_4({\RR}^{n+1})$, $G_2$ and $Spin(7)/Sp(1)$ 
(Theorem 4.1 and Corollary 5.1). But then, it
is not difficult to see that such a complex structure is not unique: since 
everything is defined by means of induced Hopf bundles,
a parallelization is induced on the fibers of a bundle in Hopf surfaces. 
For the simplest case, of $V_2({\CC}^{n+1})$, this bundle
is just the projection $V_2({\CC}^{n+1}) \rightarrow Gr_2({\CC}^{n+1})$ to the 
corresponding Grassmannian. Thus on the fibers
$S^3 \times S^1$ one can choose any complex structure that insures the 
integrability of the defined almost complex
structure on the whole total space. The family of complex structures on 
$S^3 \times S^1$ that is studied in \cite{Gau} has this
property.  

As already recalled, $V_2({\CC}^{n+1})$ admits also a family of generally 
inhomogeneous hypercomplex structures 
that contains a subfamily of homogeneneous
hypercomplex structures depending on a real parameter \cite{Bo-Ga-Ma 2}. 
As a comparison with them, we can say that all complex
structures in our family (now described by a complex parameter), project 
to the complex K\"ahler
structure of $Gr_2({\CC}^{n+1})$. Thus, since this latter complex structure 
is not compatible with the
quaternion K\"ahler structure of this Grassmannian, it follows that any of our 
complex structures on $V_2({\CC}^{n+1})$ is
non-compatible with its standard hypercomplex structures described in 
\cite{Bat 2} and \cite{Bo-Ga-Ma 2}.
\smallskip

\noindent {\bf Acknowledgement.~} We thank Paul Gauduchon for suggesting us 
to use also non-standard complex structures on our fibers $S^3 \times S^1$. 

%\bigskip\null\bigskip\null\bigskip

\section{A complex structure on some induced Hopf $S^3$-bundle}

In this paragraph we present the  key technical steps for what
follows.  
\begin{pr}\label{unu}
Let $B$ be a compact real submanifold of $\HP{n}$ and let $\pi:P\rightarrow B$ be the principal $S^3$-bundle induced over
$B$ by the  Hopf bundle $S^{4n+3}\rightarrow \HP{n}$.
If $B$ admits a Sasakian structure $(\f, \xi, \eta ,g^B)$, then one
can endow $P$ with an almost Hermitian structure. 
\end{pr}

\begin{proof}
We let $P$ have the natural pulled back metric $g$ with respect to which $\pi$
becomes a  Riemannian submersion with totally geodesic fibers
(\cite{Be}, Theorem 9.59).  
For any $X\in \mathcal{X}(B)$ we denote with $X^*$  its horizontal lift on $P$.
Let $\xi_1$, $\xi_2$, $\xi_3$ be the unit Killing vector fields which give the 
usual $3$-Sasakian structure of $S^3$ (namely, if we think about $S^3$ as
embedded in $\RR^4\cong \HH$,  $\xi_1(x)=-i x$,
$\xi_2(x)=-j x$, $\xi_3(x)=-k x$ where $i$, $j$,
$k$ are the unit imaginary quaternions) and let $\eta_1$, $\eta_2$,
$\eta_3$ be their  
duals with respect to the canonical metric of $S^3$. We regard the
$\xi_i$ 
as vector  
fields on $P$. Let $\hat{\eta}_i$ be their dual forms with respect to the 
metric $g$; their  
restrictions to any fibre coincide with the $\eta_i$. The usual 
splitting of $TP\cong \mathcal{V}\oplus\mathcal{H}$ into vertical and horizontal parts is now refined to 
$$TP \cong \text{span}\{\xi_1,\xi_2,\xi_3\}\oplus \text{span}\{\xi^*\}\oplus \mathcal{H}',$$
where $\mathcal{H}'$ represents the horizontal vector fields orthogonal to 
$\xi^*$.
We now define an almost complex structure $J$ on $P$ by:
\begin{itemize}
\item $J\xi_1=\xi_2, \quad J\xi_2=-\xi_1$,
\item $J\xi_3=\xi^*, \quad J\xi^*=-\xi_3$,
\item $JX^*=(\varphi X)^*$ for any $X\in \mathcal{X}(B)$ orthogonal to $\xi$.
\end{itemize}
As for $X\perp \xi$, $X^*$ is a section of $\mathcal{H}'$ and 
the restriction of $\varphi$ to $\xi^\perp$ is an endomorphism of
$\xi^\perp$, the last item in the definition is consistent.
One easily shows that $J^2=-1$ and is compatible with $g$. 
\end{proof}

To study the integrability of $J$ we  first compute its Nijenhuis tensor field:
$$[J,J](A_1,A_2)=[A_1,A_2]+J[JA_1,A_2]+J[A_1,JA_2]-[JA_1,JA_2],\quad
A_1,A_2\in \mathcal{X}(P).$$
As in \cite{Bo-Ga-Ma 2}, we analyse separately the different possible 
positions of $A_1$, $A_2$. We recall that, 
due to the tensorial character of $[J,J]$, 
when dealing with horizontal vector fields it is enough 
to work with basic ones whereas we always can take the $\xi_i$ as
vertical fields.

{\bf 1.} Let first $A_1=X^*$, $A_2=Y^*$. The bracket of two basic
fields $X^*,Y^*$ 
decomposes as
\begin{equation*}\label{one}
[X^*,Y^*]=[X^*,Y^*]'+\hat\eta([X^*,Y^*])\xi^*+\textrm{vertical part}.
\end{equation*}
where the $'$ denotes the $\mathcal{H}'$ part. 
By $\pi$-corelation, $[X^*,Y^*]=[X,Y]^*{}'$. Moreover, the usual
formula for the exterior derivative of a one-form 
$d\e(A_1,A_2)=A_1(\e(A_2))-A_2(\e(A_1))-\e([A_1,A_2])$
combined with $\e(X^*)=\e(Y^*)=0$ (as $X\perp \xi$ implies $X^*\perp
\xi^*$), we have
$$\e([X^*,Y^*])=-d\e(X^*,Y^*).$$
The vertical part of $[X^*,Y^*]$ must be of the form
$\sum_{i=1}^3a_i([X^*,Y^*])\xi_i$. Making the scalar product of
\eqref{one} with $\xi_j$, we find that $a_i=\e_i$. Hence
\begin{equation*}\label{3}
[X^*,Y^*]=[X^*,Y^*]'-d\e(X^*,Y^*)\xi^*-\sum d\e_i(X^*,Y^*)\xi_i.
\end{equation*}
Similarly we obtain:
\begin{equation*}
\begin{split}
[JX^*,Y^*]&=[(\varphi X)^*,Y^*]=[\varphi X,Y]^*{}'\\
&-d\e((\varphi X)^*,Y^*)\xi^*-
\sum d\e_i((\varphi X)^*,Y^*)\xi_i,\\
J[JX^*,Y^*]&=(\varphi[\varphi X,Y])^*{}'+d\e((\varphi X)^*,Y^*)\xi_3-\\
&-d\e_1((\varphi X)^*,Y^*)\xi_2+d\e_2((\varphi
X)^*,Y^*)\xi_1-d\e_3((\varphi X)^*,Y^*)\xi^*\\
J[X^*,JY^*]&=(\varphi[ X,\varphi Y])^*{}'+d\e( X^*,(\varphi Y)^*)\xi_3-\\
&-d\e_1( X^*,(\varphi Y)^*)\xi_2+d\e_2(
X^*,(\varphi Y)^*)\xi_1-d\e_3( X^*,(\varphi Y)^*)\xi^*,\\
[JX^*,JY^*]&=[\varphi X,\varphi Y]^*{}'-d\e( (\varphi X)^*,(\varphi Y)^*)\xi^*
-\sum d\e_i((\varphi X)^*,(\varphi Y)^*)\xi_i
\end{split}
 \end{equation*}
Hence we find
\begin{equation}\label{7}
\begin{split}
[J,&J](X^*,Y^*)=[\varphi X,\varphi Y]^*{}'\\
&-\left\{d\e(X^*,Y^*)-d\e((\varphi X)^*,(\varphi
Y)^*)+d\e_3((\varphi X)^*,Y^*)+d\e_3(X^*,(\varphi Y)^*)\right\}\xi^*\\
&+\left\{d\e_1((\varphi X)^*,(\varphi Y)^*)-d\e_1(X^*,Y^*)+d\e_2((\varphi
X)^*,Y^*)+d\e_2(X^*,(\varphi Y)^*)\right\}\xi_1\\
&+\left\{d\e_2((\varphi X)^*,(\varphi Y)^*)-d\e_2(X^*,Y^*)-d\e_1((\varphi
X)^*,Y^*)-d\e_2(X^*,(\varphi Y)^*)\right\}\xi_2\\
&+\left\{d\e_3((\varphi X)^*,(\varphi Y)^*)-d\e_3(X^*,Y^*)+d\e((\varphi
X)^*,Y^*)+d\e(X^*,(\varphi Y)^*)\right\}\xi_3
\end{split}
\end{equation}
As we know $[\varphi X,\varphi Y]+2d\eta(X,Y)\xi=0$ (this is the normality
condition of the Sasakian structure of $B$) the horizontal lift of
this (null) tensor field is zero, hence also its component in
$\mathcal{H}'$ is zero. But this is 
precisely $[\varphi X,\varphi Y]^*{}'$. 

On the other hand, on any Sasakian manifold one has:
$$d\eta(X,Y)=g(X,\varphi Y), \quad \varphi^2X=-X+\eta(X)\xi$$
hence $d\eta(X,\varphi Y)+d\eta(\varphi X, Y)=0$ and $d\eta(\varphi X,\varphi
Y)-d\eta( X, Y)=0$. By horizontally lifting these equations  we
get $d\e(X^*,(\varphi Y)^*)+d\e((\varphi X)^*, Y^*)=0$ and 
$d\e((\varphi X)^*,(\varphi Y)^*)-d\e(X,  Y)=0$.  Hence \eqref{7} reduces to:
\begin{equation}\label{8}
\begin{split}
[J,&J](X^*,Y^*)=\\
&-\left\{d\e_3((\varphi X)^*,Y^*)+d\e_3(X^*,(\varphi Y)^*)\right\}\xi^*\\
&+\left\{d\e_1((\varphi X)^*,(\varphi Y)^*)-d\e_1(X^*,Y^*)+d\e_2((\varphi
X)^*,Y^*)+d\e_2(X^*,(\varphi Y)^*)\right\}\xi_1\\
&+\left\{d\e_2((\varphi X)^*,(\varphi Y)^*)-d\e_2(X^*,Y^*)-d\e_1((\varphi
X)^*,Y^*)-d\e_2(X^*,(\varphi Y)^*)\right\}\xi_2\\
&+\left\{d\e_3((\varphi X)^*,(\varphi Y)^*)-d\e_3(X^*,Y^*)\right\}\xi_3
\end{split}
\end{equation}

We note that $d\e_i((\varphi X)^*,Y^*)+d\e_i(X^*,(\varphi Y)^*)=0$ iff 
 $d\e_i((\varphi X)^*,(\varphi Y)^*)-d\e_i(X^*,Y^*)=0$ (because we
can lift $\varphi$ to $P$ by defining $\hat \varphi X^*=(\varphi X)^*$ and then $\hat
\varphi$ satisfies $(\hat\varphi)^2X^*=-X^*+\e(X^*)\xi^*$).  

Hence, in order to annihilate the  $\xi^*$ and $\xi_i$ components, it
is enough to impose the condition:
\begin{equation}\label{9}
d\e_i((\varphi X)^*,(\varphi Y)^*)=d\e_i(X^*,Y^*).
\end{equation}

{\bf 2.} We now consider the case $A_1=X^*$, $A_2=\xi^*$ ($X\perp \xi$). Then 
\begin{equation*}
\begin{split}
[J,J](X^*,
\xi^*)&=[X^*,\xi^*]+J[JX^*,\xi^*]+J[X^*,J\xi^*]-[JX^*,J\xi^*]=\\
&=[X^*,\xi^*]+J[(\varphi X)^*,\xi^*]-J[X^*,\xi_3]+[(\varphi X)^*,\xi_3]
\end{split}
\end{equation*}
Here we note two wellknown facts :
\par a) On any Riemannian submersion the bracket between a
vertical field and a basic field is vertical. Hence the brackets $[X^*,
\xi_3]$ and $[(\varphi X)^*, \xi_3]$ are vertical. 
\par b) For any  connection in a  principal
bundle,  the bracket between a horizontal field and a vertical one
is horizontal. 

	As $P\rightarrow B$ is an induced $S^3$ Hopf bundle,  the horizontal 
distribution of the submersion $\mathcal{H}$ is
also the horizontal distribution of a $sp(1)$-connection (note that in
\cite{Bo-Ga-Ma 2}, when dealing with framed circle bundles, not
necessarily induced bundles, this had to be adopted as a hypothesis).
Consequently, $[X^*,
\xi_3]=[(\varphi X)^*, \xi_3]=0$. 

It remains to compute the first two terms in the expression of
$[J,J](X^*, \xi^*)$.  We have:
\begin{equation*}
\begin{split}
[X^*,\xi^*]&=[X,\xi]^*{}'-\sum d\e_i(X^*,\xi^*)\xi_i,\\
[(\varphi X)^*,\xi^*]&=[\varphi X,\xi]^*{}'-\sum d\e_i((\varphi
X)^*,\xi^*)\xi_i,\\
J[(\varphi X)^*,\xi^*]&=(\varphi [\varphi X,\xi])^*{}'\\
&-d\e_1((\varphi
X)^*,\xi^*)\xi_1+d\e_2((\varphi X)^*,\xi^*)\xi_2-d\e_3((\varphi X)^*,\xi^*)\xi_3.
\end{split}
\end{equation*}
Thus we obtain:
\begin{equation*}
\begin{split}
[J,J](X^*, &\xi^*)=([X,\xi]+\varphi[\varphi X,\xi])^*{}'
-d\e_3((\varphi X)^*,\xi^*)\xi^*\\
&+(d\e_2((\varphi X)^*,\xi^*)-d\e_1(X^*,\xi^*))\xi_1-(d\e_1((\varphi
X)^*,\xi^*)+d\e_2(X^*,\xi^*))\xi_2\\
&-d\e_3(X^*,\xi^*)\xi_3
\end{split}
\end{equation*}
We recall that on a Sasakian manifold $\varphi \xi=0$. Thus we can
add to the first paranthesis the terms $[X,\varphi \xi]-[\varphi X,\varphi \xi]$
obtaining $([X,\xi]+\varphi[\varphi X,\xi]+\varphi[X,\varphi \xi]-[\varphi X,\varphi
\xi])^*{}'$ = $([\varphi,\varphi](X,\xi))^*{}'=0$ by the normality condition
on $B$.

	Hence, in order to have $[J,J](X^*,\xi^*)=0$ it is enough to ask
\begin{equation}\label{11}
d\e_i(X^*,\xi^*)=0, \quad X\perp\xi
\end{equation}

{\bf 3.} We now choose $A_1=X^*$ and $A_2=\xi_i$ ($i=1,2$). For
$i=1$ (the other case is completely similar) we find
\begin{equation*}
\begin{split}
[J,J](X^*,
\xi_1)&=[X^*,\xi_1]+J[JX^*,\xi_1]+J[X^*,J\xi_1]-[JX^*,J\xi_1]=\\
&=[X^*,\xi_1]+J[(\varphi X)^*,\xi_1]-J[X^*,\xi_2]-[(\varphi X)^*,\xi_2]=0
\end{split}
\end{equation*}
because (see above) all the brackets are both vertical and horizontal.

	{\bf 4.} For $A_1=X^*$ and $A_2=\xi_3$ we find:
\begin{equation*}
\begin{split}
[J,J](X^*,
\xi_3)&=[X^*,\xi_3]+J[JX^*,\xi_3]+J[X^*,J\xi_3]-[JX^*,J\xi_3]=\\
&=J[X^*,\xi^*]-[(\varphi X)^*,\xi^*]
\end{split}
\end{equation*}
 The horizontal component of the reamining two brackets is 
$([\varphi[X,\xi]-[\varphi X,\xi])^*{}'-d\e_3(X^*,\xi^*)\xi_*$. 
By normality, $\eta(X)=0$,
$\varphi\xi=0$ and $d\eta(\varphi X,\xi)=0$  we have:
\begin{equation*}
\begin{split}
0=[\varphi,\varphi]&(\varphi X,\xi)+2d\eta(\varphi X,\xi)=
[\varphi X,\xi]+\varphi[\varphi^2X,\xi]+\varphi [\varphi
X,\varphi\xi]-[\varphi^2X,\varphi\xi]\\
&=[\varphi X,\xi]+\varphi[-X+\eta(X)\xi,\xi]=[\varphi[X,\xi]-[\varphi X,\xi]
\end{split}
\end{equation*}
We deduce that $([\varphi[X,\xi]-[\varphi X,\xi])^*{}'=0$, hence, as 
$d\e_3(X^*,\xi^*)=0$ according to \eqref{11}, the horizontal part of
$[J,J](X^*,\xi_3)$ is zero. Moreover, the same equation \eqref{11}
annihilates the vertical components.

{\bf 5.} Direct computation shows that in the  remaining
"mixed" case $[J,J](\xi_i,\xi^*)=0$ if 
 $[\xi_i,\xi^*]=0$
($i=1,2,3)$. As these brackets are vertical, their annulation is
equivalent with $\e_k([\xi_i,\xi^*])=0$, $k=1,2,3$. Again using the
expression of $d\e_k$ we see that we have to consider the condition:
\begin{equation}\label{12}
d\e_k(\xi_i,\xi^*)=0\quad i,k=1,2,3.
\end{equation}

{\bf 6.} We are left with the computation of $[J,J]$ on
vertical fields. Obviously $[J,J](\xi_1,\xi_2)=0$. Then
\begin{equation*}
\begin{split}
[J,J](\xi_1,\xi_3)&=[\xi_1,\xi_3]+J[J\xi_1,\xi_3]+J[\xi_1,J\xi_3]-
[J\xi_1,J\xi_3]\\
&=[\xi_1\xi3]+J[\xi_2,\xi_3]+J[\xi_1,\xi^*]-[\xi_2,\xi^*]\\
&=-2\xi_2+2J\xi_1+J[\xi_1,\xi^*]-[\xi_2,\xi^*]=0.
\end{split}
\end{equation*}
by \eqref{12}. The case $A_1=\xi_2$, $A_2=\xi_3$ is similar.
Summing up we have proved:
\begin{pr}\label{pp}
The following conditions are sufficient for the  almost complex structure 
defined in Proposition
\ref{unu}  to be integrable:

	1)  $d\e_k(\xi_i,\xi^*)=0\quad i,k=1,2,3.$

	2)  $d\e_i(X^*,\xi^*)=0$, for any $X\perp \xi$ and $i=1,2,3$.

	3) $d\e_i((\varphi X)^*,(\varphi Y)^*)=d\e_i(X^*,Y^*)$ for any $X,	Y\perp \xi$ and $i=1,2,3$.
\end{pr}

Observe now that $d\e_k$ can be
identified
as the vertical parts of the curvature form $\Omega$ of
the $sp(1)$ connection $\mathcal{H}$. Moreover:
\begin{pr}
$\mathcal{H}$ is an $sp(1)$ connection if and only if the vector
fields $\xi_i$ are Killing on $(P,g)$.
\end{pr}

\begin{proof}
Recall that $\mathcal{H}$ is a connection iff for any $X\in
\Gamma(\mathcal{H})$ and any vertical $V$, the brackets $[X,V]$ are
horizontal.  As any horizontal field is of the form $aX^*+b\xi^*$, we
have $[a\xi^*+bX^*,V]=a[\xi^*,V]-V(a)\xi^*+b[X^*,V]-V(b)X^*$ hence
$[a\xi^*+bX^*,V]$ is horizontal iff $[\xi^*,V]$ and $[X^*,V]$ are
horizontal. We can take $V=\xi_i$. The above two brackets are surely
vertical, thus they will be horizontal iff they are zero. 

Let us compute the Lie derivative of the
metric $g$ on the total space in the direction $\xi_i$. We obtain
successively:
$$(\Ll_{\xi_i}g)(X^*,\xi^*)=\xi_ig(X^*,\xi^*)-g([\xi_i,X^*],\xi^*)-g(X^*,[\xi_i,\xi^*]=0$$
because $g(X^*,\xi^*)=0$ and the brackets in the last two terms are
vertical.
$$(\Ll_{\xi_i}g)(X^*,Y^*)=\xi_ig(X^*,Y^*)-g([\xi_i,X^*],Y^*)-g(X^*,[\xi_i,Y^*]=0$$
as $g(X^*,Y^*)$ does not depend on vertical directions and again
because the brackets in the last two terms are
vertical.
$$(\Ll_{\xi_i}g)(X^*,\xi_k)=\xi_ig(X^*,\xi_k)-g([\xi_i,X^*],\xi_k)-g(X^*,[\xi_i,\xi_k]$$
Here $g(X^*,\xi_k)=0$, $[\xi_i,\xi_k]2\epsilon_{ikj}\xi_j$ and
$g(X^*,\xi_j)=0$. Hence 
$$(\Ll_{\xi_i}g)(X^*,\xi_k)=-g([\xi_i,X^*],\xi_k)=-\e_k([\xi_i,X^*])=d\e_k(\xi_i,X^*).$$
$$(\Ll_{\xi_i}g)(\xi^*,\xi_k)=-g([\xi_i,\xi^*],\xi_k)=-\e_k([\xi_i,\xi^*])=d\e_k(\xi_i,\xi^*).$$
We obtained that $\xi_i$ are Killing fields iff $[\xi_i,X^*]$ and
$[\xi_1,\xi^*]$ are horizontal.
\end{proof} 
>From the proof we also obtained that condition 1) of the above
proposition is assured. We can finally give the integrability
condition of the constructed $J$ in terms of curvature properties of
$\mathcal{H}$. 

\begin{te}\label{intcur}
The  almost complex
structure in Proposition \ref{unu} is integrable if the curvature form
of the $sp(1)$ connection $\mathcal{H}$ satisfies the conditions:
\begin{itemize}
\item $\Omega((\varphi X)^*,(\varphi Y)^*)=\Omega(X^*,Y^*)$ for any $X,	Y\perp \xi$ and $i=1,2,3$.
\item $\Omega(X^*,\xi^*)=0$, for any $X\perp \xi$ and $i=1,2,3$.
\end{itemize}
\end{te}

We may observe that the stated conditions express the compatibility
between the Sasakian structure of the base (which is not induced by
the immersion of $B$ in $\mathbb{H}P^n$) and the bundle structure of
$P\rightarrow B$.

\begin{re}	

(i)  The K\"ahler form 
$\omega$ of $(P,g,J)$ is 
 non-closed, and indeed it does not satisfy any of the Gray-Hervella 
conditions besides the
integrability of $J$. A similar computation proves that 
$L_{\xi^*}J=L_{\xi_3}J=0$, thus $\xi^*$ and $\xi_3$ are
infinitesimal automorphisms of the constructed complex structure. 

(ii) We note also that by its definition the complex structure $J$ on
$P$ depends on the choice of a the $3$-Sasakian structure of 
$S^3$. Different choices of the $3$-Sasakian 
triples $\{\xi_1,\xi_2,\xi_3\}$ define complex structures that are conjugated 
in $End(TP)$. More informations about the dependence of $J$ on the
chosen parallelization of $S^3$ are given in \S 4 for the case of
$V_2(\CC^{n+1})$ and $\tilde V_4(\RR^{n+1})$.

(iii) Although the construction of $J$ does not use explicitely the
induced Hopf bundle, the construction doesn't work for merely 
Riemannian submersions with fibres $S^3$: one needs a canonical way of choosing the parallelization of $S^3$ in order to
avoid monodromy problems.

\end{re}

%\bigskip 

\section{The zero level sets of two moment maps}

Consider now the two maps 
$$\mu:\HH^{n+1}\rightarrow Im~\HH, \qquad \nu:\HH^{n+1}\rightarrow Im~Ê\HH^3,$$
defined in the coordinates $h=[h_0:h_1:...:h_n]$ of $\HH^{n+1}$ by 
$$\mu(h)=\sum_{a=0}^n\ov{h}_aih_a, \quad \quad
\nu(h)=(\sum_{a=0}^n\ov{h}_aih_a, \sum_{a=0}^n\ov{h}_ajh_a,\sum_{a=0}^n\ov{h}_akh_a),$$
and recall that $\mu$ and $\nu$ can be interpreted as the moment maps associated to the diagonal action of
$U(1)$ and of $Sp(1)$ on the 3-Sasakian sphere $S^{4n+3} \subset \HH^{n+1}$ fibering over $\HH P^{n}$. The corresponding
quaternion K\"ahler reductions  are the quaternion K\"ahler Wolf spaces $SU(n+1)/S(U(n-1)\times U(2))\cong
Gr_2(\CC^{n+1})$ and
$SO(n+1)/(SO(n-3)\times SO(4))\cong\widetilde{Gr}_4(\RR^{n+1})$, respectively (cf. for example \cite{Bo-Ga 2}). 

We proved in \cite{orp} the following:
\begin{pr}\label{uu}
 \par {\rm (i)} 
$\mu^{-1}(0)$ is diffeomorphic to 
 the total space of the induced Hopf $S^1$-bundle via the
Pl\"ucker embedding $Gr_2(\CC^{n+1}) \hookrightarrow  \CC P^N$, and consequently a Sasakian metric
is induced on $\mu^{-1}(0)$ by the Pl\"ucker embedding of this Grassmannian. 
\par {\rm (ii)} The zero level set $\nu^{-1}(0)$ is diffeomorphic to the total space
of the induced Hopf $S^1$-bundle over the Fano manifold
$Z_{\widetilde{Gr}_4(\RR^{n+1})}$, by means of the 
embeddings $Z_{\widetilde{Gr}_4(\RR^{n+1})}\hookrightarrow Gr_2(\CC^{n+1})\hookrightarrow 
\CC P^N$, the first of which is defined by regarding $Z_{\widetilde{Gr}_4(\RR^{n+1})}$ as the space of totally isotropic
two-planes in $\CC^{n+1}$.  Thus an induced Sasakian metric is obtained on $\nu^{-1}(0)$.  
\end{pr}

Since both $\mu^{-1}(0)$ and $\nu^{-1}(0)$ can be shown to be simply
connected, the first statement both of (i) and of (ii) is a
consequence of the following observation: {\it Let $\pi:P\rightarrow
B$ be a principal circle bundle with simply connected $P$ over a smooth complex algebraic projective submanifold $B$ of
$\CC P^N$ with
$H^2(B,\ZZ)\cong \ZZ$. 
Then $P$ is diffeomorphic to 
the total space of the induced Hopf bundle of $S^{2N+1}\rightarrow
\CC P^N$, via the embedding $B \hookrightarrow  \CC P^N$.}
In the case of $\mu^{-1}(0)$, the submanifold $B$ is the Grassmannian $Gr_2(\CC^{n+1})$ and its Pl\"ucker embedding is
used in $\CC P^N$, $N=\binom{n+1}{2}-1$.
As for $\nu^{-1}(0)$, it is also an induced Hopf $S^1$-bundle but over the twistor space 
$Z_{\widetilde{Gr}_4(\mathbb{R}^{n+1})}$ of the quaternion K\"ahler real Grassmannian
$\widetilde{Gr}_4(\mathbb{R}^{n+1})$. This twistor space is a complex submanifold of
$Gr_2(\CC^{n+1})$ \cite{Kobak}. On the other hand, the  composition of the fiberings
$$\nu^{-1}(0)\stackrel{S^1}{\rightarrow} Z_{\widetilde{Gr}_4(\mathbb{R}^{n+1})}
\stackrel{S^2}{\rightarrow}\widetilde{Gr}_4(\mathbb{R}^{n+1})$$  is a
$SO(3)$-bundle which endows $\nu^{-1}(0)$ with a $3$-Sasakian structure
via the inversion theorem 4.6 of \cite{Bo-Ga 1}.
%\bigskip

\section{Applications to Stiefel manifolds}

If we regard the Stiefel manifolds $V_2(\CC^{n+1})$ and $\widetilde
V_4({\RR}^{n+1})$ as homogeneous manifolds, we immediately recognize
them as total spaces of the induced  bundles $S^3\rightarrow
S^{4n+3}\rightarrow \HH P^{n}$ over $\mu^{-1}(0)$, respectively. The conditions 
stated in Theorem \ref{intcur} are verified for these bundles (cf. \cite{orp}). 
This gives the following:

\begin{te}
A family of uncountably many homogeneous complex structures on the Stiefel manifolds $V_2(\CC^{n+1})$ and $\widetilde
V_4({\RR}^{n+1})$ can be obtained by combining the K\"ahler-Einstein structures of $Gr_2(\CC^{n+1})$ and 
$Z_{\widetilde{Gr}_4(\mathbb{R}^{n+1})}$ with any of the complex structure on the Hopf surface $\CC^2 -\{0\}/(z
\rightarrow
\lambda z$), given by all choices of $\lambda \in \CC^*,  \vert \lambda \vert >1$.
\end{te}

\begin{proof} A standard complex structure on $V_2(\CC^{n+1})$ and $\widetilde
V_4({\RR}^{n+1})$ is obtained by applying Proposition \ref{unu} and Theorem \ref{intcur} 
to the highest vertical arrows in the diagram:
\begin{equation*}
\begin{array}{cccccl}
\widetilde{V}_4(\mathbb{R}^{n+1})&\;  \hookrightarrow\;&
V_2(\mathbb{C}^{n+1})&\; \hookrightarrow\;&S^{4n+3}&\\[1mm]
\Big\downarrow  {\scriptstyle S^3}& &\Big\downarrow  {\scriptstyle S^3}& &\Big\downarrow  {\scriptstyle S^3} \\[2mm]
\nu^{-1}(0)&\; \hookrightarrow\;&\mu^{-1}(0)&
\hookrightarrow\;&\mathbb{H}P^n&\;\;\;\; S^{2N+1}\\[1mm]
\Big\downarrow  {\scriptstyle S^1}& &\Big\downarrow  {\scriptstyle S^1}& &&\!\!\!\!\!\swarrow \\[2mm]
Z_{\widetilde{Gr}_4(\mathbb{R}^{n+1})}&\;
\hookrightarrow\;&Gr_2(\mathbb{C}^{n+1})&\;
\hookrightarrow\;&\mathbb{C}P^N&\\[2mm]
\Big\downarrow  {\scriptstyle S^2}&&&&\;\\[2mm]
\widetilde{Gr}_4(\mathbb{R}^{n+1})&&&&&
\end{array}
\end{equation*}
where Proposition 3.1 is applied to zero level sets $\mu^{-1}(0)$ and
$\nu^{-1}(0)$ to obtain their induced Sasakian structures on them.

The same diagram tells that $V_2(\CC^{n+1})$ and $\widetilde
V_4({\RR}^{n+1})$ are bundles in Hopf surfaces $S^3 \times S^1$ over
the complex K\"ahler-Einstein manifolds 
$Gr_2(\mathbb{C}^{n+1})$ and $Z_{\widetilde{Gr}_4(\mathbb{R}^{n+1})}$, respectively. 
On all these fibers $S^3 \times S^1$ a 
simultaneous parallelization is induced by a choice of a 
3-Sasakian structure on $S^{4n+3}$ and a Sasakian structure on
$S^{2N+1}$. From this point of view, the above mentioned complex 
structure on the Stiefel manifold is by construction given by the
choice of the standard complex structure on the fibers $S^3 \times S^1$ and by 
the lift of the complex structure of the
K\"ahler-Einstein bases.  The integrability of the whole complex
structure was 
insured by the computations described above.
%, and
%by noting that, by regarding now our Stiefel manifolds as
%$(S^3\times S^1)$-bundles, the only non trivial part of these computations is 
%the vanishing of the mixed components $[J,J](\xi^*,X^*)$
%and $[J,J](\xi_\alpha,X^*)$ of the Nijenhuis tensor. 
%To check that these components are zero, we used the
%Sasakian structures involved. 

Observe now  that these same computations, leading to $[J,J]=0$,  
can be carried out even if the complex structure on the fibers is not 
defined in the standard
way (i.e. $J\xi_1=\xi_2, \quad J\xi_2=-\xi_1, \quad J\xi_3=\xi^*,
\quad J\xi^*=-\xi_3$), but according to formulas like:
\begin{align*}
J\xi_1=\xi_2, &\quad J\xi_2=-\xi_1 \\ 
J\xi^*= \alpha \xi^* + \beta \xi_3, & \quad J\xi_3= \gamma \xi^* + \delta \xi_3,
\end{align*}
where the matrix $\pmatrix \alpha, \beta\\ 
\gamma, \delta 
\endpmatrix$, whose entries are real and constant, has trace $0$ and
determinant $1$. The complex 
structures defined in this way on $S^3 \times S^1= \CC^{2} -
\{0\}/(z \rightarrow \lambda z)$ correspond to all the possible 
choices of the generator $\lambda \in
\CC^* = \CC -\{0\}, \vert \lambda \vert >1$, and it can be shown that 
all these complex structures on the Hopf surface are
inequivalent (cf. \cite{Gau}, p. 142-143).
\end{proof}

 Note that these complex structures on $V_2(\CC^{n+1})$ project to the complex 
structure with
respect to which the symmetric Grassmannian $Gr_2(\CC^{n+1})$ is
K\"ahlerian. But this Grassmannian also has a quaternion K\"ahler
structure which does not contain the K\"ahler
structure (\emph{i.e.} whilst the K\"ahler metric coincides with the
quaternion-K\"ahler one, the complex structure compatible with the
K\"ahler metric is not a section of the quaternion bundle). 
On the other hand, it is the quaternion K\"ahler structure of 
the Grassmannian $Gr_2(\CC^{n+1})$ that
produces, \emph{via} the associated homogeneous 3-Sasakian 
manifold and its deformations, the hypercomplex structures on
$V_2(\CC^{n+1})$
\cite{Bat 2},
\cite{Bo-Ga-Ma 2}. This gives the following:

\begin{co}
The constructed complex structures on $V_2(\CC^{n+1})$ are
non-com\-pa\-ti\-ble with its standard hypercomplex structure. 
\end{co}
%\bigskip\null\bigskip

\section{Two special cases}

More complex structures on Stiefel manifolds can be obtained by looking at the 
following exceptional cases.
Observe that the group $G_2$ can be regarded as the "special" 
Stiefel manifold of 
coassociative orthonormal 4-frames $(e_1,e_2,e_3,e_4)$ in
$\RR^7$. This means that the corresponding 4-plane has an orthogonal 
complement that is an associative 3-plane, \emph{i.e.} closed under
the vector product of $\RR^7$. This follows easily from the references 
\cite{Ma}, p. 252, \cite{Ha-La}, p.115. The second
reference states in fact that $G_2 \cong V_3^\phi (\RR^7)$, 
the latter being the Stiefel manifold of orthonormal 3-frames
$(e_1,e_2,e_4)$ such that, with respect to the product of Cayley numbers, 
$e_4 \perp e_1 e_2$. Of course such 3-frames are in
one-to-one correspondence with coassociative 4-frames via 
$(e_1,e_2,e_4) \leftrightarrow (e_1,e_2,e_1 e_2, e_4)$. The Stiefel
manifold $G_2$ fibers in Hopf surfaces $S^3 \times S^1$ over the flag manifold 
$G_2/U(2)^+$, twistor space of the quaternion
K\"ahler submanifold
$G_2/SO(4)$ of
$\widetilde{Gr}_4(\RR^7)$.
\par Also related to the geometry of Cayley numbers, 
one can consider the "special" Stiefel manifold of
Cayley 4-frames in $\RR^8$, \emph{i.e.} orthonormal 4-frames  
spanning 4-planes in $\RR^8$ that are closed under the double
cross-product (cf. again \cite{Ma}, p. 
261, \cite{Ha-La}, p. 118). The Stiefel manifold of Cayley 4-frames is easily 
recognized to be the homogeneous space
$Spin(7)/Sp(1)$, fibering again in Hopf surfaces over the twistor space of the 
Grassmannian of Cayley 4-planes
$Spin(7)/(Sp(1)\times Sp(1)
\times Sp(1))/\ZZ_2)$. This latter manifold is a quaternion K\"ahler 
submanifold of $\widetilde{Gr}_4(\RR^8)$.
\par This discussion extends to the homogeneous 3-Sasakian bundles and 
yields the following two diagrams of
submanifolds considered in more detail in  \cite{orp2}. The first diagram is:
\begin{equation*}
\begin{array}{ccccccl}
V=G_2&\; \hookrightarrow\;&\widetilde{V}_4(\mathbb{R}^{7})&\;  \hookrightarrow\;&
V_2(\mathbb{C}^{7})&\; \hookrightarrow\;&S^{27}\subset \HH^7 \\[1mm]
\Big\downarrow  {\scriptstyle S^3}& &\Big\downarrow  {\scriptstyle S^3}& &\Big\downarrow  {\scriptstyle S^3}&
&\Big\downarrow  {\scriptstyle S^3} \\[2mm] G_2/Sp(1)^+&\; \hookrightarrow\;&\nu^{-1}(0)&\; \hookrightarrow\;&\mu^{-1}(0)&
\hookrightarrow\;&\mathbb{H}P^6\\[1mm]
\Big\downarrow  {\scriptstyle S^1}& &\Big\downarrow {\scriptstyle S^1}& &\Big\downarrow {\scriptstyle S^1}& &\\[2mm]
G_2/U(2)^+&\; \hookrightarrow\;&Z_{\widetilde{Gr}_4(\mathbb{R}^{7})}&\;
\hookrightarrow\;&Gr_2(\mathbb{C}^{7})&&\\[2mm]
\Big\downarrow {\scriptstyle S^2}& &\Big\downarrow  {\scriptstyle S^2}&&&&\;\\[2mm]
G_2/SO(4)&\; \hookrightarrow\;&\widetilde{Gr}_4(\mathbb{R}^{7}),&&&&
\end{array}
\end{equation*}
\noindent where the + sign appearing in the left column corresponds to a choice that is 
significant for the structure of the two homogeneous
manifolds $G_2/Sp(1)^+$ and $G_2/U(2)^+$, cf. \cite{Sa}, p. 164.

Similarly, one gets a second diagram by considering Cayley 4-frames and 
Cayley 4-planes in $\RR^8$:
\begin{equation*}
\begin{array}{ccccccl}
V=\frac{Spin(7)}{Sp(1)}&\; \hookrightarrow\;&\widetilde{V}_4(\mathbb{R}^{8})&\;  \hookrightarrow\;&
V_2(\mathbb{C}^{8})&\; \hookrightarrow\;&S^{31}\subset \HH^8 \\[1mm]
\Big\downarrow  {\scriptstyle S^3}& &\Big\downarrow  {\scriptstyle S^3}& &\Big\downarrow  {\scriptstyle S^3}&
&\Big\downarrow 
{\scriptstyle S^3} \\[2mm]
\frac{Spin(7)}{Sp(1)\times Sp(1)}&\; \hookrightarrow\;&\nu^{-1}(0)&\; \hookrightarrow\;&\mu^{-1}(0)&
\hookrightarrow\;&\mathbb{H}P^7\\[1mm]
\Big\downarrow  {\scriptstyle S^1}& &\Big\downarrow  {\scriptstyle S^1}& &\Big\downarrow  {\scriptstyle S^1}& &\\[2mm]
\frac{Spin(7)}{(Sp(1)\times Sp(1)\times U(1))/\ZZ_2}&\; \hookrightarrow\;&Z_{\widetilde{Gr}_4(\mathbb{R}^{8})}&\;
\hookrightarrow\;&Gr_2(\mathbb{C}^{8})&&\\[2mm]
\Big\downarrow  {\scriptstyle S^2}& &\Big\downarrow  {\scriptstyle S^2}&&&&\;\\[2mm]
\frac{Spin(7)}{(Sp(1)\times Sp(1)\times Sp(1))/\ZZ_2}&\; \hookrightarrow\;&\widetilde{Gr}_4(\mathbb{R}^{8}).&&&&
\end{array}
\end{equation*}
\bigskip

These two diagrams, combined with Proposition \ref{unu} and 
Theorem \ref{intcur}, give:
\begin{co} An uncountable family of homogeneous complex structures is obtained on the special Stiefel manifolds $G_2$ and
$Spin(7)/Sp(1)$, by regarding them as induced Hopf bundles of $S^{27} \rightarrow \HH P^6$ and of $S^{31}Ê\rightarrow \HH
P^7$ over the Sasakian submanifolds $G_2/Sp(1)^+ \subset \HH P^6$, 
$\frac{Spin(7)}{Sp(1)\times Sp(1)} \subset \HH P^7$, respectively.
\end{co}

\end{document}